\documentclass[a4paper, 10pt]{article}
\usepackage{graphicx}
\usepackage{enumitem}
\usepackage{amsfonts}
\usepackage{stmaryrd}
\usepackage{textcomp}

\begin{document}
\begin{center}
\textbf{LINKS BETWEEN BROWNIAN PROCESSES, PROPORTIONAL INCREMENTS AND THE LAST HITTING TIMES}
\vspace*{0.5\baselineskip}\\
Meziane Privat\\
Université Libre de Bruxelles
\end{center}
\normalsize
\vspace*{1.5\baselineskip}
\begin{footnotesize}
\textbf{Abstract :} In this short article, we will focus on the different links between some stochastic processes resulting from Brownian motion and two notions of probability theory (proportional increments and last hitting times).
\vspace*{0.5\baselineskip}
\vspace*{0.5\baselineskip}
\end{footnotesize}
\vspace*{1.5\baselineskip}\\
{\large \textbf{1.~~Introduction}}
\vspace*{1\baselineskip}\\
In this work, we will first focus on the link between stochastic processes resulting from a Brownian motion, such as Brownian motion with drift and proportional increments processes introduced by Bruss and Yor in 2012 in [BY12], where these authors used the new notion of proportional increments processes to solve an optimal stopping time problem. We will demonstrate the converse of a theorem they state in [BY12]. \\
In the second part of this work, we will focus mainly on another article of Bruss and Yor [BY15] published in 2015. We will follow their work on the last hitting times, they used to demonstrate the Williams decomposition of the dimension 3 Bessel process ($ BES (3) $).\\ 
We will end with a reflection on the interest of these mathematical researches in the context of compassionate use clinical trials.
\vspace*{1\baselineskip}\\
{\large \textbf{2.~~Brownian motion and proportional increments}}
\vspace*{1\baselineskip}\\
\textbf{Definition 2.1 :} A Brownian motion from $ X_0 $ is a process of the form $ X_t = X_0 + B_t $ (with $ B_t $ a standard Brownian motion), with as initial condition $ X_0 $ independent of $ (B_t)_{t\in\mathbb{R}}$.
\vspace*{0.3\baselineskip}\\
\textbf{Definition 2.2 :} A Brownian motion derived from $ X_0 $ with drift (or trend) $ a $ and diffusion coefficient $ \sigma $ is a process of the form $ X_t = X_0 + \sigma B_t + at $ with the initial condition  $ X_0 $ independent of $ (B_t)_{t \in \mathbb {R}}$. \\
\begin{center}
\includegraphics[scale=0.3]{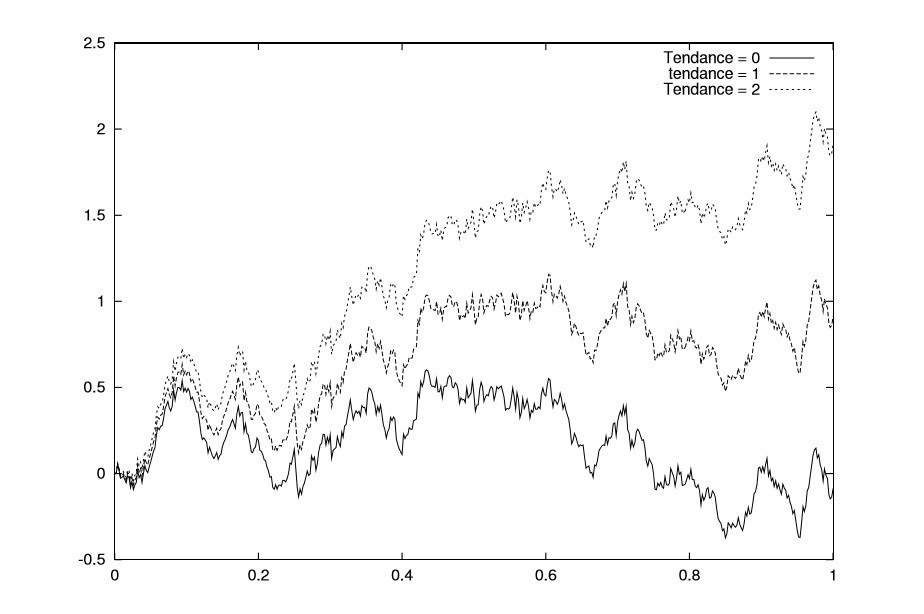} \\
Figure 1 : Representation of three path of a Brownian motion with drift for $X_0=0$ and $\sigma =1$
\end{center}
\begin{flushleft}
{\footnotesize Source : ($https://www.editions.polytechnique.fr/files/pdf/EXT\_1579\_4.pdf$)}
\end{flushleft}
\textbf{Definition 2.3 :} Let $ N_t $ be a stochastic process defined on a filtered probability space $ (\Omega, \mathcal {F}, (\mathcal {F} _t), \mathbb {P}) $ with natural filtration $ \mathcal {F} _t = \lbrace N_u : u \leq t \rbrace$. \\ We say that $ (N_t) _ {t> 0} $ is a process with proportional increments on $] 0, + \infty [$ if,
\begin{center}
$\forall t>0$ with $N_t\neq 0, \forall s\geq 0 : E(N_{t+s}-N_t\vert\mathcal{F}_t)=\frac{s}{t}N_t~ a.s.$
\end{center}
\textbf{Remark :} Let $ (X_t) _ {t \in \mathbb {R}} $ a Brownian motion with drift $ a \neq 0 $, then $ (X_t) _ {t \in \mathbb {R}} $ is not a process with proportional increments.
\vspace*{0.3\baselineskip}\\
\textbf{Proof (of the remark) :} Let $ (X_t) _ {t \in \mathbb {R}} $ a Brownian motion with drift $ a \neq 0 $ and diffusion coefficient $ \sigma $ ($ X_t = X_0 + \sigma B_t + at $).\\
Let $ t $ such that $ X_t \neq 0 $ and let $ s \geq 0 $,\\
\begin{center}
$E(X_{t+s}-X_t\vert\mathcal{F}_t)=E(X_0+\sigma B_{t+s}+a(t+s)-(X_0+\sigma B_t+at)\vert\mathcal{F}_t)~a.s.$\\
$=E(\sigma B_{t+s}-\sigma B_t+as\vert\mathcal{F}_t)~a.s.$~~~~~~~~~\\
$=E(\sigma (B_{t+s}-B_t)\vert\mathcal{F}_t)+as~a.s.$~~~~~~~~~\\
$=E(\sigma (B_{t+s}-B_t))+as~a.s.$,~~~~~~~~~~~~\\ by independence of the increments.\\
$~~~~~~=as~a.s.$ since $E(\sigma (B_{t+s}-B_t))=0$~~~~~~~~~\\
$\neq\frac{s}{t}X_t~a.s.~~~~~~~~~~~~~~~~~~~~~~~~~~~~~~~~~~~$
\end{center}
\begin{flushright}
$\square$
\end{flushright}
\textbf{Notation :} We will introduce for the rest of this part the process $ (X_t^*) _ {t \in \mathbb {R}} $ defined by $ \forall t \in \mathbb {R} ~ X_t ^ * = X_0 ^ * + \sigma tB_t + at $ with $ B_t $ the standard Brownian motion and $ a, \sigma \in \mathbb {R}$.
\vspace*{0.3\baselineskip}\\
\textbf{Property 2.1 :} The process $ (X_t ^ *) _ {t \in \mathbb {R}}$ is a proportional increment process if $ X_0 ^ * = 0$.
\vspace*{0.3\baselineskip}\\
\textbf{Proof :} Let $ t $ such that $ X ^ * _ t \neq 0 $ and let $ s \geq 0 $ and $X_0^*=0$,\\
\begin{center}
$E(X_{t+s}^*-X_t^*\vert\mathcal{F}_t)=E(X_0^*+\sigma (t+s)B_{t+s}+a(t+s)-(X_0^*+\sigma tB_t+at)\vert\mathcal{F}_t)~a.s.$\\
$=E(\sigma (t+s)B_{t+s}-\sigma tB_t+as\vert\mathcal{F}_t)~a.s.$~~~~~~\\
$~~~~~~~=E(\sigma t(B_{t+s}-B_t)+\sigma sB_{t+s}\vert\mathcal{F}_t)+as~a.s.$~~~~~~~~~\\
$~~~=E(\sigma (B_{t+s}-B_t))+E(\sigma sB_{t+s}\vert\mathcal{F}_t)+as~a.s.$,\\ by linéarity of the expectation.\\
$~~~~~~~~~~~~~~~~~~~~=E(\sigma sB_{t+s}\vert\mathcal{F}_t)+as~a.s.$, since $E(\sigma t(B_{t+s}-B_t))=0$~~~~~~~\\
$~~~~~~~~~~~=\sigma sB_{t}+as~a.s.$, since $(B_{t})_{t>0}$ is a martingale~~~~~~~~\\
$=\frac{s}{t}(0+\sigma tB_t+at)~a.s.$~~~~~~~~~~~~~~~~~~~~~~~~~~~\\
$=\frac{s}{t}X_t^*~a.s.$ since  $X_0=0$~~~~~~~~~~~~~~~~~~~~~~~~~~
\end{center}
So $\forall t>0$ with $X^*_t\neq 0, \forall s\geq 0 : E(X^*_{t+s}-X^*_t\vert\mathcal{F}_t)=\frac{s}{t}X^*_t~ a.s.$
\begin{flushright}
$\square$
\end{flushright}
\textbf{Theorem 2.1 :} Let $ (N_t) _ {t> 0} $ be a proportionally incremental process and let $ R_t = \frac {N_t} {t} $. If $ E (\vert R_t \vert) <\infty, (R_t)_{t>0}$ is a martingale. \\
This theorem is demonstrated in [BY12]. We will show that the reciprocal is true.\vspace*{0.3\baselineskip}\\
\textbf{Reciprocal 2.1 :} Let $ (M_t) _ {t> 0} $ be a martingale such that $ E (\vert tM_t \vert) <+\infty $, then the process $ N_t = tM_t $ is a proportionally incremental process.
\vspace*{0.3\baselineskip}\\
\textbf{Proof :} Let $ (M_t) _ {t> 0} $ be a martingale. Let $N_t=tM_t$.\\
\begin{center}
Alors $\forall t\neq 0, M_t\neq 0~a.s.\Rightarrow N_t\neq 0~a.s.$
\end{center}
Let $ t> 0 $ such that $ N_t \neq 0 $ and let $ s \geq 0 $ then,
\begin{center}
$E(N_{t+s}-N_t\vert\mathcal{F}_t)=E((t+s)M_{t+s}-tM_t\vert\mathcal{F}_t)~a.s.$\\
$~~~~~~~~~~~~~~~~~~~~~~~~~~~~=E(t(M_{t+s}-M_t)+sM_{t+s}\vert\mathcal{F}_t)~a.s.$\\
$~~~~~~~~~~~~~~~~~~~~~~~~~~~~~~~~~~~~~~=E(t(M_{t+s}-M_t)\vert\mathcal{F}_t)+E(sM_{t+s}\vert\mathcal{F}_t)~a.s.$,\\ ~~~~~~~~~~~~~~~~~~~~~~~~~~~~~~~~by linéarity of the expectation\\
So $E(N_{t+s}-N_t\vert\mathcal{F}_t)~~= sM_t$ since $(M_t)_{t>0}$ is a martingale\\
$=\frac{s}{t}tM_t~a.s.$~~~~\\
$=\frac{s}{t}N_t~a.s.~~~~~$
\end{center}
\begin{flushright}
$\square$
\end{flushright}
\textbf{Property 2.2 :} The process $ (tB_t) _ {t> 0} $ with $(B_t)_{t>0}$ the standard Brownian motion, is a proportionally incremental process.
\vspace*{0.3\baselineskip}\\
\textbf{Proof :} The standard Brownian motion is a martingale and $ E (\vert tB_t \vert) <+\infty $ so according to the previous reciprocal the property is proved.
\begin{flushright}
$\square$
\end{flushright} 
\textbf{Property 2.3 :} Let $ (X_t ^ {* (k)})_{t>0}$ with $ k \in \llbracket 1, n \rrbracket :=\lbrace 1,2,...,n\rbrace$ a sequence of stochastic processes defined as in Property 2.1. Suppose $\forall~t>0~X_t ^ {* (k)}$ is $\mathcal{F}_t$-measurable. Put $ S_t = \sum \limits ^ n_ {k = 1} c_kX_t ^ {* (k)} $ with $ \forall k \in \llbracket 1, n \rrbracket ~ c_k $ arbitrary constants. Then $ (S_t) _ {t> 0} $ is a proportional increment process.
\vspace*{0.3\baselineskip}\\
\textbf{Proof :} $ S_t $ is $ \mathcal {F} _t $-measurable as a sum of $ \mathcal {F} _t $-measurable process. \\
In addition, the condition $ (S_t \neq 0) $ is met. \\
Let's calculate $ E (S_ {t + s} -S_t \vert \mathcal {F} _t)$,
\begin{center}
$E(S_{t+s}-S_t\vert\mathcal{F}_t)=E(\sum\limits^n_{k=1} c_kX_{t+s}^{*(k)}-\sum\limits^n_{k=1} c_kX_t^{*(k)}\vert\mathcal{F}_t)~a.s.$\\
$~~~~~~~~~~~~~~=\sum\limits^n_{k=1} c_kE(X_{t+s}^{*(k)}-X_t^{*(k)}\vert\mathcal{F}_t)~a.s.$\\
$=\sum\limits^n_{k=1} c_k\frac{s}{t} X_{t}^{*(k)}~a.s.$~~~~ \\ ~~~~~~~~~~~~~~~~~~~~~~~~~~~~~~~~~~~~~~since $\forall k\in\llbracket 1,n\rrbracket ~E(X_{t+s}^{*(k)}-X_t^{*(k)}\vert\mathcal{F}_t)=\frac{s}{t}X_{t}^{*(k)}~a.s.$\\
$=\frac{s}{t}S_t~a.s.$~~~~~~~~~~~~~~~
\end{center}
\begin{flushright}
$\square$
\end{flushright}
\textbf{Remark :} This property is extended to all proportional incremental processes by Bruss and Yor in [BY12].
\vspace{1\baselineskip}\\
{\large \textbf{3.~~ 
Williams' decomposition of the three-dimensional Bessel process with a look at the latest hitting time}}
\vspace{1\baselineskip}\\
\normalsize
For the following presentation we will follow mainly the article of Bruss and Yor of 2015 [BY15].
\vspace{0.1\baselineskip}\\
\textbf{Motivation :} Determine if a stochastic process hits a certain set for the last time depends on what will happen in the future. Therefore, the latest hitting times are generally not measurable with respect to the natural process filtration and hence among the difficult random moments of a stochastic process. Downtime, on the other hand, has this property by definition, and we have a rather impressive collection of theorems and tools for downtime. As Chung concludes (see the quote from Nikeghbali and Platen (2013) [NP13]), it is imperative to avoid the latest hitting times.\\
It's a way of looking at things, but quite often the reality is a little different. Bruss and Yor argue that many interesting problems related to the theory of optimal shutdown require us to take into account the latest typing times, not the downtime. The attitude has therefore changed and the works of Jeulin (1980) [Jeu80] and others have strongly influenced this development. In our article you will find interesting examples from the field of mathematical finance, and we will explore one of them from a different angle. We would also like to broaden slightly the horizon of mathematical finance by looking at some other examples.
\vspace{0.5\baselineskip}\\
\textbf{3.1~~The $ BES (3) $ process and Williams' theorem}
\vspace{0.5\baselineskip}\\
\textbf{Definition 3.1 :} The Bessel process of order $ n $ is the process $(X_t)_{t>0} $ given by,
\begin{center}
$X_t = \Vert B_ {t}\Vert$
\end{center}
where $ \Vert . \Vert $ denotes the Euclidean norm in $ \mathbb {R} ^ n $, and $ (B_t) _ {t> 0} $ is a Brownian motion with $ n $ dimensions starting at the origin.
\vspace{0.3\baselineskip}\\
\textbf{Definition 3.2 :} The Bessel process at $ n $ dimensions is the solution of the stochastic differential equation,
\begin{center}
$dX_ {t} = dB_ {t} + {\frac {n-1}{2}} {\frac{dt}{X_{t}}}$
\end{center}
where $(B_t)_{t\geq 0}$ is a Brownian motion with one dimensions.\\
We note $BES(n)$.
\vspace{0.3\baselineskip}\\
\begin{center}
\includegraphics[scale=0.1]{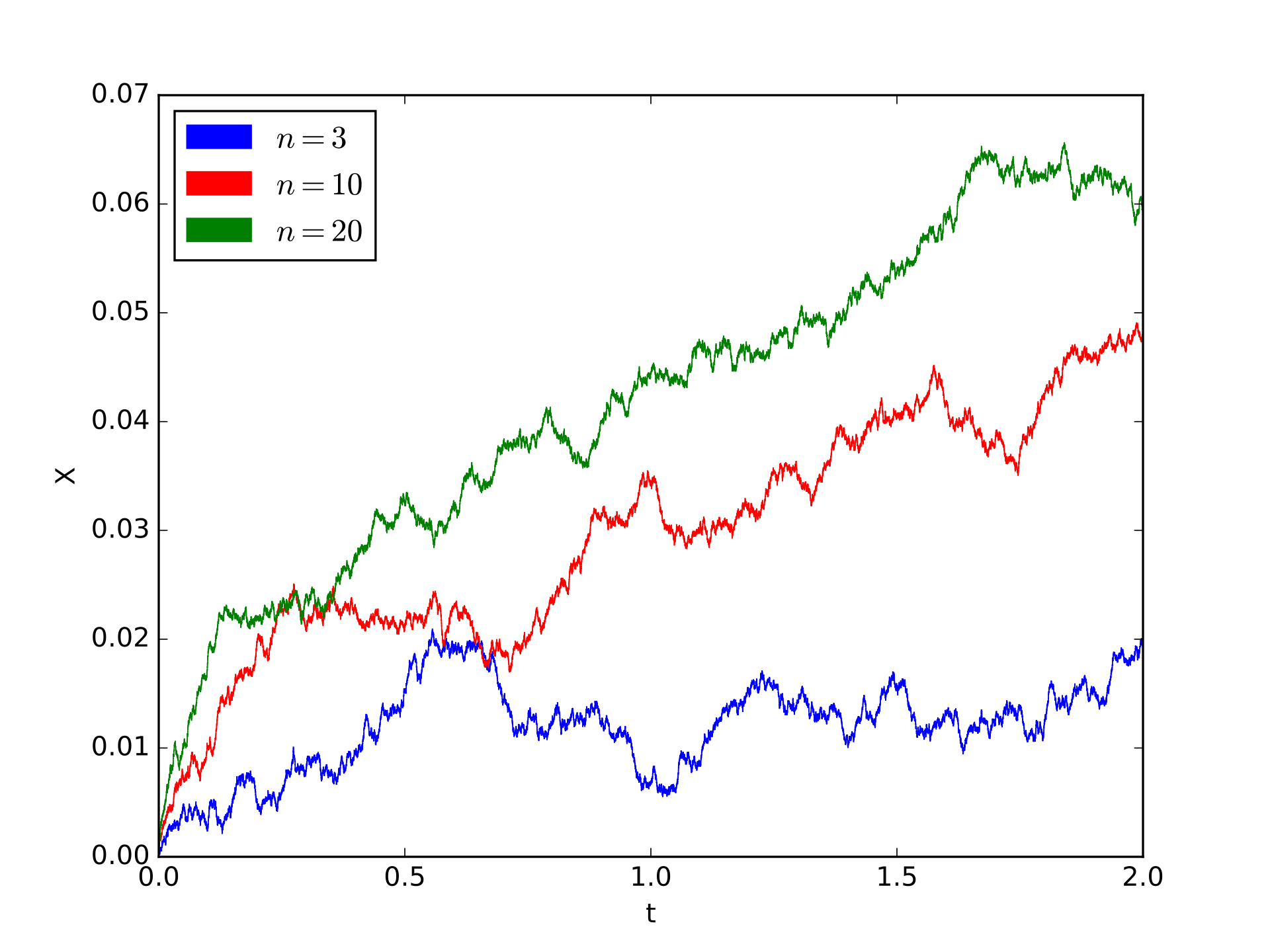} \\
Figure 2 : Three achievements of the Bessel processes.
\end{center}
\begin{flushleft}
{\footnotesize (Source : $https://en.wikipedia.org/wiki/Bessel\_process$)}
\end{flushleft}
Let $ (R_t) _ {t \geq 0} $ be a $3$-dimensional Bessel process ($ BES (3) $ process) on $ \mathbb {R} _ + $ from $ r> 0 $. Note by $ (\mathcal {F} _t) _ {t \geq 0} $ its natural filtration and let $ I_t $ be the infimum of the process $(R_t)_{t>0} $ at time $ t $, that is to say,
\begin{center}
$I_t=\inf\limits_{s\leq t} R_t$
\end{center}
The following results can be found in Nikeghbali and Platen (2013) [NP13] around corollary 4.10, \\
\begin{itemize}
\item [{(i)}]: $ I_ \infty $ follows the same law as the random variable $ rU $, where $ U $ is a uniform law over $ [0,1] $.
\item [{(ii)}]: The Azema's supermartingale associated with the random time $ g $ at which the process $ (R_t)_{t>0} $ reaches $ I_ \infty $ is given by,
\begin{center}
$ Z_t \equiv P (g> t \vert \mathcal {F} _t) = \frac {I_t} {R_t} $.
\end{center}
\textbf{Reminder :} We call Azema's martingale the process $(\mu_t)_{t\geq 0}$ such that, \begin{center}
$\forall t \geq 0, \mu_t = E[B_t\vert\mathcal{F}_t]$
\end{center} and Azema's supermartingale if, \begin{center}
$\forall t \geq 0, \mu_t \geq E[B_t\vert\mathcal{F}_t]$.
\end{center}With $ (B_t) _ {t> 0} $ a standard Brownian motion.
\item [{(iii)}]: The Laplace transform of the $ g $ law is,
\begin{center}
$ E (e ^ {- \lambda g}) = \frac {1} {\sqrt {2 \lambda} r} (1-e ^ {\sqrt {2 \lambda} r}) $
\end{center}
\item [{(iv)}]: The density of $ g $ is given by $ p (t) $ whom is equal to,
\begin{center}
$p(t)=\frac{1}{\sqrt{2\pi t}r}(1-e^{\frac{r^2}{2t}})$.
\end{center}
\end{itemize}
The objective of Bruss and Yor (2015) [BY15] is now to show Williams' decomposition of a $ BES (3) $ process to its ultimate minimum, and how this decomposition is closely related to (i) - (ii ) - (iii) - (iv).
\vspace{0.5\baselineskip}\\
Recall that if $ (B_t) _ {t \geq 0} $ is a Brownian motion starting from $ 0 $ and $ a $ is a constant, then the law of the first hitting time of $ a $ by $(B_t)_{t>0}$ , denoted $ T_a ^ {(B)} $, is given by,
\begin{center}
\hfill $P(T^{(B)}_a\in dt) = \frac{dt}{\sqrt{2\pi t^3}}\vert a\vert e^{-\frac{a^2}{2t}}$\hfill(1). 
\end{center}
This well-known fact allows us to rewrite statements (iii) and (iv) above as follows,
\begin{center}
\hfill$g =^{\mathcal{L}} T^{(B)}_{rU}$\hfill(2)
\end{center}
where $ U $ is independent of $(B_t)_{t>0}$ and uniform over $ [0,1] $, and where $ = ^ {\mathcal {L}} $ is equality in law. This can be verified with (iii) and (iv). In fact, (2) can be understood through the classical decomposition of the process $(R_t)_{t>0}$ before and after the time $ g $, thanks to Williams (1974) [Wil74]. 
This is explained later.
\vspace{1\baselineskip}\\
\textbf{3.2~~Williams decomposition of $(R_t)_{t>0}$, before and after $ g $, by progressive enlargement}
\begin{center}
\includegraphics[scale=0.9]{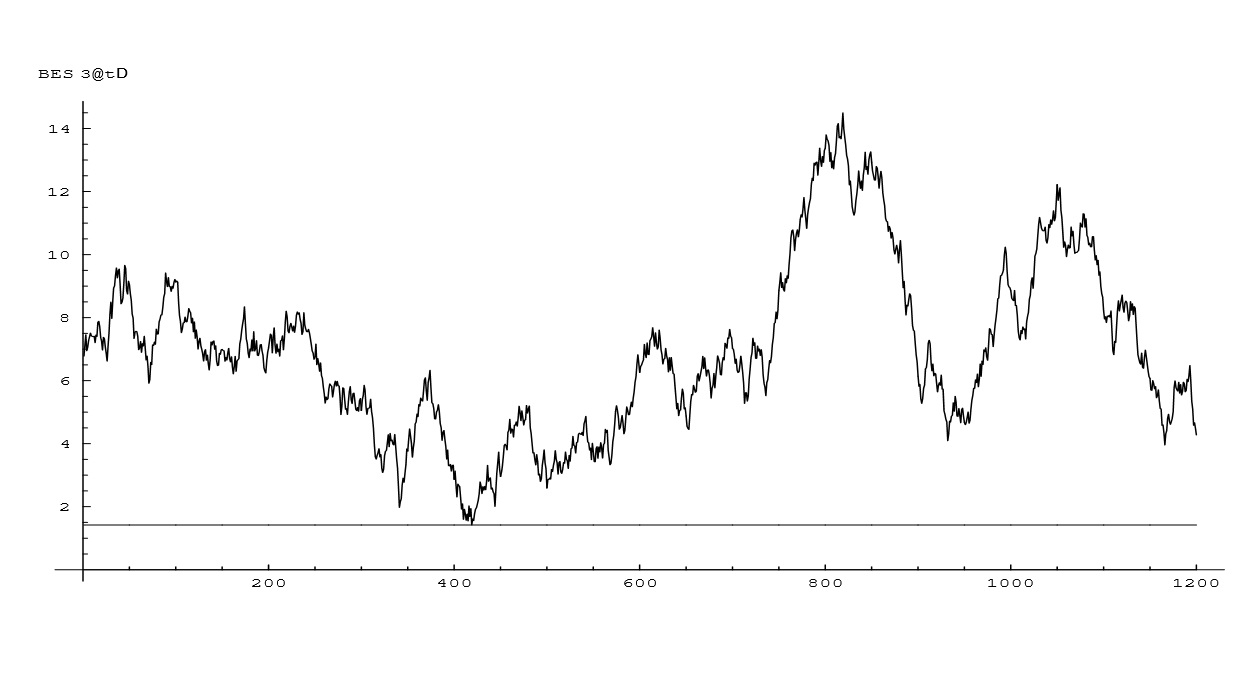} \\
Figure 3 : Decomposition of a $ BES (3) $ process 
\end{center}
Note that this figure is nothing other than a version (simulated infinite horizon) of Figure 5 in Revuz-Yor (1999) [RY99] (see Proposition 3.10 and Theorem 3.11 of Chapter 6, Section 3) where the process $ BES (3) $ is considered from the level $ c:= r $.
\vspace{0.5\baselineskip}\\
We now state precisely Williams' decomposition theorem before and after the $ g $ time.
\vspace{0.3\baselineskip}\\
\textbf{Theorem 3.1} (Williams (1974) [Wil74])\\
Consider the following three  random variables, that we suppose independent,
\begin{itemize}
\item [{(i)}] a Brownian motion $ (B_t ') _ {t \geq 0} $ with $ B_0' = r> 0 $;
\item [{(ii)}] a uniform random variable $ U $ on $ [0,1] $;
\item [{(iii)}] a process $ BES (3) $ $ (\tilde {R} _t) _ {t \geq 0} $ with $ \tilde {R} _0 = 0 $;
\end{itemize}
Let $ (R) $ be the process defined by,
\begin{center}
\hfill$R_t$=
$\left\lbrace
\begin{array}{c c}
B_t',$~~~~~~~~~~~~~if~$t\leq g\\ 
rU+\tilde{R}_{t-g}$~~~~ else $
\end{array}
\right.$ \hfill(3)
\end{center}
with $ g = \inf \lbrace u \geq 0: B_u '= rU \rbrace $. So $(R_t)_{t>0}$ is a $ BES (3) $ process starting with $ r> 0 $. \\
We note that the pre-$g$-Browien motion found in (3) explains the result (2). Indeed, if $ B_t '= r-B ^ {(0)} _ t $ with,
\begin{center}
\hfill$ g = \inf\lbrace u\geq 0: B^{(0)}_u = r(1-U)\rbrace$,\hfill (2') 
\end{center} 
so (2') implies (2).
\vspace{0.5\baselineskip}\\
Bruss and Yor (2015) [BY15] then move to the proof of the theorem via the magnification formula which describes the additive decomposition of the process $ BES (3) (R_t) $ in the filtration $ (\mathcal {F} _t ^ g) $ with filtration $ (\mathcal {F} _t) $ and taking $ g $ as the stopping time. \\
First, we have,
\begin{center}
\hfill$R_t = r + B_t + \int^t_0 \frac{ds}{R_s},$\hfill(4)
\end{center}
where $ (B_t) $ is a Brownian motion with respect to $ (\mathcal {F} _t) $. \\
Secondly, the enlargement formula (see for example Jeulin (1980) [Jeu80]) gives,
\begin{center}
\hfill$r + B_t = B_t' + \int_0^{g\wedge t}\frac{d <B, Z>_u}{Z_u}+ \int_g^t\frac{d <B, 1-Z>_u}{1-Z_u}$\hfill(5)
\end{center}
with $ (B'_t) $ a Brownian motion relative to $ (\mathcal {F} ^ g_t) $. \\
Third, we deduce from (ii) the two following identities,
\begin{center}
\hfill$\frac{d <B, Z>_u}{Z_u}=-\frac{du}{R_u}$ si $u\leq g$\hfill(6)\\
et\hfill~\\
\hfill$\frac{d <B, 1-Z>_u}{1-Z_u}=\frac{I_\infty du}{R_u(R_u-I_\infty)}$ si $u>g$\hfill(7)
\end{center} 
These two identities involve (using (4) and (6)), as well as 
\begin{center}
$\frac{1}{R_u}+\frac{I_\infty}{R_u(R_u-I_\infty)}=\frac{1}{R_u-I_\infty}$
\end{center}
the form of the process pre- $g$ and the form of the process post-$ g $. \\
Finally, for the proof of (3) to be complete, they prove that the process $(B_t')_{t>0}$ is independent of the random variable $ I_\infty = ^ \mathcal {L} rU $, or more precisely since $ I_ \infty = a $, the process pre-$g$ is just the process $ (B'_u) _ {u \leq T'_a} $ with an obvious notation.
\vspace{1\baselineskip}\\
\textbf{{\large 4~~More general context of last hitting times}}
\vspace{1\baselineskip}\\
Now, as announced in the introduction, we are going to focus on compassionate use, but before that we will present different contexts related to the last hitting times brought by Bruss and Yor (2015) [BY15].
\vspace{0.3\baselineskip}\\
Buying and selling problems, choice problems, secretary problems and others are typical representatives of a last-minute problem. They can be seen as stopping problems on the last "improvement" of a stochastic process. In some of these problems, the difficulty with the character of the last hit disappears. To give a very simple example, suppose that we observe sequentially the variables $ X_1, X_2, ... $ and want to maximize them, for a given objective function $ f $, the expected total return, so we're looking for,
\begin{center}
$\arg\max\limits_\tau f(X_1,X_2,...,X_\tau)$.
\end{center}
Suppose now that the optimal payment for stopping after time $ t $ does not depend on $ \mathcal {F} _t $, where $ (\mathcal {F} _s) $ indicates natural filtration. So,
\begin{center}
$\sup\limits_{\tau\geq t} E(f(X_1,X_2,...,X_\tau\vert\mathcal{F}_t)) = \sup\limits_{\tau\geq t} E(f(X_{t+1},X_{t+2},..., X_\tau)).$ 
\end{center}
so that the RHS and $ X_1, X_2, ..., X_t $ are $ \mathcal {F} _t $ -measurable. Thus, we have to compare at each instant $ t $ the value of $ f (X_1, X_2, ..., X_t) $ with the supremum of the RHS in order to make the optimal decision.
In more difficult problems, independence with $ \mathcal {F} _t $ is usually no longer satisfied. However, external information about the underlying process can sometimes help to solve the problem of the last hitting time into a tractable problem. For example, we can refer to the proof of the law $ \frac{1}{e} $ of the best choice of Bruss (1984) [Bru84] which transformed the sequential problem of the best choice into a combinatorial problem in a context non-sequential and there are other examples that we could see as simple because well defined as we can find in [Bru00], [Den13] or [Tam10]. \\
The other extreme is the example of the Robbins problem where the complexity is quite different because the optimal strategy depends at all times on all the information. For this problem we can refer to section 4.2 of Bruss and Ferguson in [BF96]. \\
In addition, there are some other problems where the last typing time can be interesting to use, such as discovering the first time a random subset of a given set is complete. They give a concrete example. This is in the important field of clinical trials, specifically in compassionate clinical trials.
\vspace{0.5\baselineskip}\\
\textbf{4.1~~ Compassionate use clinical trials}
\vspace{0.5\baselineskip}\\
In such trials, a group of patients is treated with a drug not yet on the market or whose doses have not been validated by the authorities, and which may have significant side effects, the only justification is usually that this may be the last hope for patients.\\
In general, little is known about the likelihood of success of unapproved drugs or unapproved doses of (known) drugs. For this reason, compassionate use trials are usually sequential so that physicians or statisticians can learn from previous observations. It should also be mentioned that such tests may not be allowed in all countries, even with the patient's approval.\\
All the reflections show that the successive treatments cause cultural problems. Conscientious physicians must try to save all the lives that can be saved and, at the same time, avoid unnecessary suffering caused by treatment. Since the doctor is not a prophet, the goal should be to stop (in a given group of patients in a given horizon) with maximum probability, the first patient ending the random subset of success that stops with the last success. Indeed, all successes are covered, while the remaining patients (de facto not recoverable by the drug) must not suffer unnecessarily. Usually, the treatment sequence is also interrupted if the current estimate of the probability of success falls below a certain limit. No detailed proposal is made regarding the estimation method. In [Bru18] several case studies are proposed.\\
If the doctor has incomplete information on the respective probabilities of success, the general solution of the optimal stopping problem is an open problem and, as the authors think, an important problem.
However, it is interesting to note that if the physician has no information, the optimal solution can be deduced from the notion of proportional increments introduced in Bruss and Yor (2012) [BY12], as shown by Dendievel (2013) [Den13] in a closely related problem.\\
Proposals and mathematical methods should be of interest to physicians because it would streamline their behavior towards patients they can no longer save using traditional methods. Not to mention that when a life is at stake, the doctor must do everything to optimize its decisions. \\
Applying this kind of method would also minimize the number of unnecessary future treatments by stopping at the optimal time. \\
Moreover, one can wonder whether these different mathematical methods can not be introduced into the artificial intelligences that will shape medicine in the coming years. Whether to help the doctor by indicating different possibilities of action or by telling him what to do.
\vspace{1\baselineskip}\\
\noindent
{\large \textbf{Acknowledgement}}
\vspace{1\baselineskip}\\
\normalsize
I am grateful to Professor Franz Thomas Bruss, from whom I discovered the interest of proportional increments processes and of the last hitting times problems.
\vspace{2\baselineskip}\\
{\large \textbf{References}}
\normalsize
\begin{enumerate}
\item [{[AKM10]}] K. Ano, H. Kakinuma, et N. Miyoshi, 2010, \textit{Odds theorem with multiple selection chances},J. Appl. Prob., Vol. 47, n° 4, pages 1093-1104.
\item [{[Bru84]}] F.T. Bruss, 1984, \textit{A unified approach to a class of best-choice problems with an unknown number of options}, Annals of Probab., Vol. 12, n°3, pages 882-889. 
\item [{[Bru00]}] F.T. Bruss, 2000, \textit{Sum the odds to one and stop}, Annals of Probab., Vol. 28, n°3, pages 1384-1391.
\item [{[Bru18]}] F.T. Bruss, 2018, \textit{A Mathematical Approach to Comply with Ethical Constraints in Compassionate Use Treatments}, Mathematical Scientist, Vol. 43, n°1, pages 10-22.
\item [{[BF96]}] F.T. Bruss et T.S. Ferguson, 1996, \textit{Half-Prophets and Robbins Problem of Minimizing the Expected Rank}, Athens Conference on Applied Probability and Time Series Analysis, Vol 1, Applied Probability, Lecture Notes in Statistics 114, Springer, pages 1-17.
\item [{[BY12]}] F.T. Bruss et M. Yor, 2012, \textit{Stochastic processes with proportional increments and the last-arrival problem}, Stochastic Process and their Application n°122, pages 3239-3261, Elsevier.
\item [{[BY15]}] F.T. Bruss, et M. Yor, 2015, \textit{A New Proof of Williams’ Decomposition of the Bessel Process of Dimension Three with a Look at Last-Hitting Times}, Bulletin of the Belgian Mathematical Society - Simon Stevin, Belgian Mathematical Society, n°22 (2), pages 319-330.
\item [{[Den13]}] R. Dendievel, 2013, \textit{New developments of the odds-theorem}, The Math. Scientist, Vol.38, n°2, pages 111-123.
\item [{[Fer16]}] T.S. Ferguson, 2016, \textit{The Sum-the-Odds Theorem with Application to a Stopping Game of Sakaguchi}, Mathematica Applicanda, n°44 (1), pages 45-61.
\item [{[Jeu80]}] T. Jeulin, 1980, \textit{Semi-martingales et grossissement d’une filtration}, Lect.Notes Math. 833, Springer Verlag. 
\item [{[NP13]}] A. Nikeghbali et E. Platen, 2013, \textit{A reading guide for last passage times with financial applications in view}. Finance and Stochastics, Vol. 17, n°3, pages 615-640.
\item [{[RY99]}] D. Revuz and M. Yor, 1999, \textit{Continuous martingales and Brownian motion}. 3rd Edition, Springer.
\item [{[Tam10]}] M. Tamaki, 2010, \textit{Sum the multiplicative odds to one and stop}, J. Appl. Prob., Vol. 47, n°3, pages 761-777.
\item [{[Wil74]}] D. Williams, 1974, \textit{Path decomposition and continuity of local times for one-dimensional diffusions, I}, Proceed. of London Math. Soc. (3), 28, pages 738-768.
\end{enumerate}
\end{document}